\documentclass{article}

\usepackage{arxiv}

\usepackage[utf8]{inputenc} 
\usepackage[T1]{fontenc}    
\usepackage{hyperref}       
\usepackage{url}            
\usepackage{booktabs}       
\usepackage{amsfonts}       
\usepackage{nicefrac}       
\usepackage{microtype}      
\usepackage{cleveref}       
\usepackage{graphicx}
\usepackage[numbers]{natbib}
\usepackage{doi}

\newtheorem{theorem}{Theorem}
\newtheorem{Proposition}[theorem]{Proposition}%

\newtheorem{Example}{Example}%

\title{Classes of second order nonlinear partial differential equations reducible to first order}

\author{Noureddine Mhadhbi \\
	Department of Mathematics\\
	 College of Sciences and Arts\\
	King Abdulaziz University \\ P.O.Box 344, Rabigh Campus,21911, Saudi Arabia \\
	\texttt{nalmhadhbi@kau.edu.sa} \\
	\And
	Sameh Gana\\
	Department of Basic Sciences\\
	Deanship of Preparatory Year and Supporting Studies\\
Imam Abdulrahman Bin Faisal University\\P.O. Box 1982, Dammam,34212, Saudi Arabia\\
	\texttt{sbgana@iau.edu.sa} \\
		\And
	Hamad Khalid Alharbi\\
	Mathematics Department, Faculty of Science, King Abdulaziz University\\ P. O. Box 80203,  Jeddah 21589, Saudi Arabia.\\
	\texttt{Halharbi0502@stu.kau.edu.sa}\\
Applied College\\
	Imam Mohammad Ibn Saud Islamic University, Saudi Arabia\\
	\texttt{hakhalharbi@imamu.edu.sa} \\
}

\begin{document}
\maketitle

\begin{abstract}
	In this paper, we present new techniques for solving a large variety of partial differential equations. The proposed method reduces the PDEs to first order differential equations known as classical equations such as Bernoulli, Ricatti and Abel equations. The main idea is based on implementing new techniques by combining variations of parameters with characteristic methods to obtain many new and general exact solutions.  In each class of PDE's, we give illustrated examples. Moreover, the method presented in this paper can be easily extended to classes of second order nonlinear PDEs.
\end{abstract}


\section{Introduction}

Partial differential equations are widely used as models to describe complex physical phenomena in various areas of science, particularly in fluid mechanics, solid-state physics and plasma waves. Nonlinear differential equations cover also optics, bio-hydrodynamics and nonlinear quantum mechanics. \\ 
In recent years, we notice significant progress in the expansion of effective methods for obtaining exact solutions of 
nonlinear partial differential equations, we can mention the inverse scattering transform ~\citep{Ablowitz}, the extended tanh-function method ~\citep{Fan}, the truncated expansion method ~\citep{Kudryashov}, F-expansion method ~\citep{D.Zhang,Zhang}, Jacobi elliptic method ~\citep{Shikuo},  the Backlund transformations ~\citep{Hong}, the sine-cosine function method ~\citep{Wazwaz}, the (G'/G)-expansion method ~\citep{Wang} and the extended Kudryashov method ~\citep{Hassan}.

One of the effective tools to solve PDEs is the method of characteristics ~\citep{Polyanin2}, ~\citep {Brian}, ~\citep{Myint}, ~\citep{Zachmanoglou} and ~\citep {Rhee}. The kernel point of this method is the idea that the partial differential equation and the given initial conditions can be reduced into a set of differential equations which can be directly and easily solved and we obtain general exact solutions of the initial problem.

The method of variation of parameters ~\citep{Mahouton}, ~\citep{Polyanin1}, ~\citep{Jovan}, ~\citep{Olver} and ~\citep{Kevorkian} can also be successfully used in some cases to reduce by one the order of nonlinear
partial differential equations. In general, reducing the order of partial differential equations enables us to find suitable analytical methods for the resolution.   

However, some classes of nonlinear second order partial differential equations can not be exactly solved using classical methods. As a consequence, we need to recur to novel methods for solving a large variety of PDEs.

In this paper, we start by presenting the results which are relatively known in the simplest case, the solutions of linear PDEs in section 2 and we establish the exact solutions of nonlinear PDEs in section 3. The ideas will extend to more complicated cases and several examples are presented for illustration.\\
In sections 4 and 5, we solve some classes of PDEs called Bernoulli type and Ricatti Model:

$$
	u_{t}+a(x,t)u_{x}=b(x,t)u+\alpha(x,t)u^{n},
$$
and 
$$
	u_{t}+a(x,t)u_{x}=b(x,t)u+\alpha(x,t)+\beta(x,t)u^{2}. 
$$
 In sections 6,7 and 8, we prove that the variation of parameters is successfully implemented in combination with the characteristic method to obtain the exact solutions of the nonlinear second order 
 partial differential equation of the form: 
$$u_{tt}+a(x,t)u_{xt}=b(x,t)+(u_{t}+a(x,t)u_{x})f(u), $$
$$u_{xt}+a(x,t)u_{xx}=b(x,t)+(a(x,t)u_{x}+u_{t})f(u),$$
$$f'(u_{t})(u_{tt}+au_{xt})=B(x,t)+A(u) (u_{t}+au_{x}).$$
In the last section, we extend the proposed method to general classes of nonlinear second order 
partial differential equation of the form: 
$$u_{tt}+a(x,t)u_{xt}+b(u)u_{t}(u_{t}+a(x,t)u_{x})=\alpha(x,t)e^{-\int b(u)du} +G(u) (u_{t}+a(x,t)u_{x}).$$
The special cases $b(u)=-\frac{1}{u}$ and  $G(u)=\beta u^{n}$ have been investigated in ~\cite{Mahouton}. We apply our approach for finding more general exact solutions.

\section{Partial Differential Equations of the form 
		$u_{t}+a(x,t)u_{x}-\alpha(x,t)u=b(x,t)$ }
	
We consider the first order partial differential equation of the following form:
\begin{equation}
	\left\{
	\begin{array}{c}
		u_{t}+a(x,t)u_{x}-\alpha(x,t)u=b(x,t) \\
		u(x,0)=\phi(x).
	\end{array}
	\right.
\end{equation}
where $u$ is a function of $(x,t)\in \mathbb{R}^{2}$.
\\
First, we solve the differential equation $\left\{
\begin{array}{c}
	\frac{dx(t)}{dt}=a(x(t),t) \\
	x(0)=x_{0}.
\end{array}
\right. $ \\
Then (1) can be rewritten as \begin{equation}u_{t}(x(t),t)+a(x(t),t)u_{x}-\alpha(x(t),t)u=b(x(t),t).
\end{equation}
We multiply (2) by $e^{-\int_{0}^{t}\alpha(x(s),s)ds}$, we get 
\begin{equation}\frac{d}{dt}[e^{-\int_{0}^{t}\alpha(x(s),s)ds}u]=b(x(t),t)e^{-\int_{0}^{t}\alpha(x(s),s)ds}.
\end{equation}
Therefore, the following statement holds:
\begin{Proposition}
	The first order partial differential equation (1) can be reduced to the differential equation (3) and the solution of (1) is
	$$u(x,t)=\phi(x_{0})e^{\int_{0}^{t}\alpha(x(s),s)ds}+e^{\int_{0}^{t}\alpha(x(s),s)ds}\int_{0}^{t}b(x(\tau),\tau)e^{-\int_{0}^{\tau}\alpha(x(s),s)ds}d\tau.$$
\end{Proposition}
For illustration, let us consider this example:
\begin{Example} Let $a(x,t)=x$, $\alpha(x,t)=1$ and $b(x,t)=x+t$.
	\begin{equation}
		\left\{
		\begin{array}{c}
			u_{t}+xu_{x}-u=x+t \\
			u(x,0)=\phi(x).
		\end{array}
		\right. 
	\end{equation}
\end{Example}
Let $x(t)$ be the solution of $\left\{
\begin{array}{c}
	\frac{dx(t)}{dt}=x \\
	x(0)=x_{0}.
\end{array}
\right. $

then $x(t)=x_{0}e^{t}$.

(4) can be rewritten as 
$$\frac{d}{dt}u(x(t),t)=u(x(t),t)+x+t.$$
Then $$\frac{du}{dt}-u=x+t,$$
$$e^{-t}\frac{du}{dt}-e^{-t}u=(x(t)+t)e^{-t},$$
$$\frac{d}{dt}(e^{-t}u)=(x(t)+t)e^{-t}.$$
Therefore
\begin{center}
	\begin{tabular}{c c c }
		$u(x(t),t)$& = &$e^{t}\phi(x_{0})+e^{t}\int_{0}^{t}(x(s)+s)e^{-s}ds$ \\
		& = &$e^{t}\phi(xe^{-t})+e^{t}\int_{0}^{t}(x_{0}e^{s}+s)e^{-s}ds$.
	\end{tabular}
\end{center}
The solution of (4) is $$u(x,t)=e^{t}\phi(xe^{-t})+xt-t-1+e^{t}.$$
\begin{figure}[hbt!]
	\centering
	\includegraphics[width=7cm]{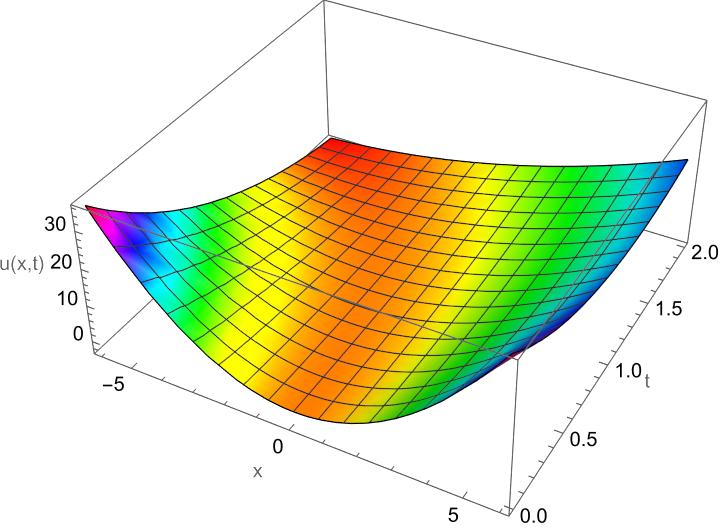} 
	\caption{Solution of (4) 
		with $\phi(x)=x^{2}$. \label{fig2}}
\end{figure}

\section{Partial Differential Equations of the form $
		u_{t}+a(x,t)u_{x}=f(u)b(x,t)$ }
	
\vspace{0.3cm}
Let us consider the first order partial differential equation of the following type:
\begin{equation}
	\left\{
	\begin{array}{c}
		u_{t}+a(x,t)u_{x}=f(u)b(x,t) \\
		u(x,0)=\phi(x).
	\end{array}
	\right.
\end{equation}
Let $x(t)$ be the solution of $\left\{
\begin{array}{c}
	\frac{dx(t)}{dt}=a(x(t),t) \\
	x(0)=x_{0}.
\end{array}
\right. $ \\
Then (5) is transformed to  $$\frac{d}{dt}u(x(t),t)=f(u(x(t),t))b(x(t),t),$$
and \begin{equation}\frac{du}{f(u)}=b(x(t),t)dt. 
\end{equation}
Then we deduce the following result:
\begin{Proposition}
	The first order partial differential equation (5) can be reduced to equation (6). Furthermore, let $F(u)=B(t)$ 
	where $B(t)=b(x(t),t)$ and 
	$F(u(x_{0},0))=b(x_{0},0)$, then the general solution of (5) is easily determined by solving the equation $$F(u)=F(\phi(x_{0}))+\int_{0}^{t}b(x(s),s)ds.$$
\end{Proposition}
As an illustration, let us consider the following examples:
\\
\begin{Example}Let $a(x,t)=x$, $f(u)=u^{2}$ and $b(x,t)=1$.
	\begin{equation}
		\left\{
		\begin{array}{c}
			u_{t}+xu_{x}=u^{2} \\
			u(x,0)=x.
		\end{array}
		\right.  
	\end{equation}
\end{Example}
Since $\frac{du}{u^{2}}=1$,
we get $$u(x(t),t)=\frac{1}{-t+\frac{1}{xe^{-t}}}.$$

\begin{figure}[hbt!]
	
	\centering
	\includegraphics[width=7cm]{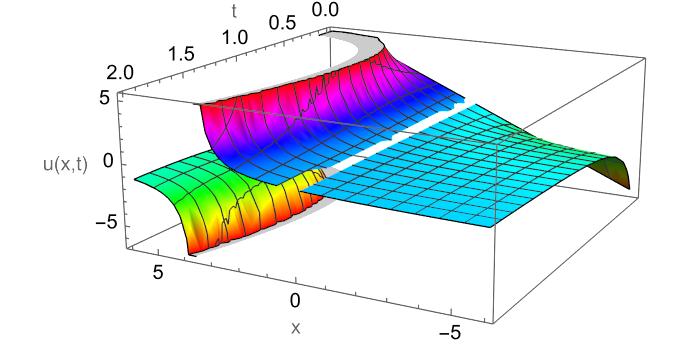}  
	
	\caption{Solution of (7) 
		with $\phi(x)=x$.\label{fig2}}
\end{figure}
\begin{Example} Let $a(x,t)=x$, $f(u)=u^{2}+1$ and $b(x,t)=x+t^{2}$.
	\begin{equation}
		\left\{
		\begin{array}{c}
			u_{t}+xu_{x}=(u^{2}+1)(x+t^{2}) \\
			u(x,0)=\phi(x).
		\end{array}
		\right.  
	\end{equation}
\end{Example}
We get $$\frac{du}{dt}=(u^{2}+1)(x+t^{2}),$$
$$\tan^{-1}(u)=\tan^{-1}(\phi(xe^{-t}))+\int_{0}^{t}(x_{0}e^{s}+s^{2})ds.$$
Then $$u(x,t)=\tan(\tan^{-1}(\phi(xe^{-t}))+x(1-e^{-t})+\frac{t^{3}}{3}).$$
\begin{figure}[h]
	\centering 
	\includegraphics[width=7cm]{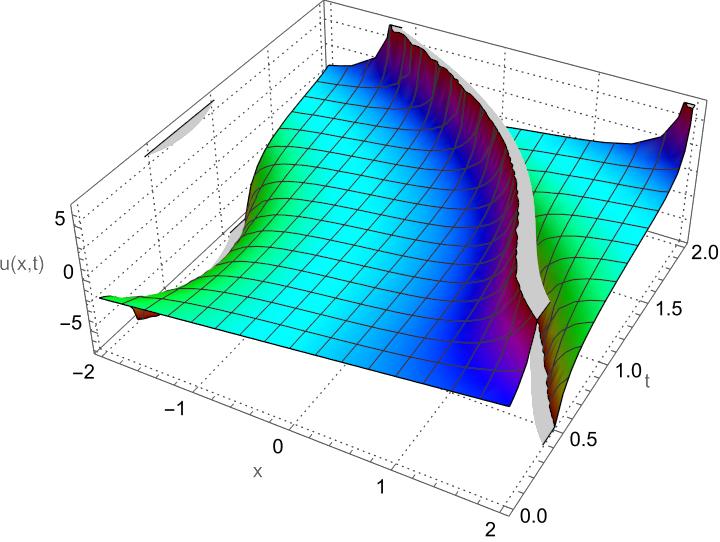}
	\caption{Solution of (8) with $\phi(x)=x.$\label{fig2} }	
\end{figure}
\section{Bernoulli Equation of the form $
		u_{t}+a(x,t)u_{x}=b(x,t)u+\alpha(x,t)u^{n}$}
	
\vspace{0.3cm}
We consider the Bernoulli Equation of the form:
\begin{equation}
	\left\{
	\begin{array}{c}
		u_{t}+a(x,t)u_{x}=b(x,t)u+\alpha(x,t)u^{n}\\
		u(x,0)=\phi(x).
	\end{array}
	\right.
\end{equation}
Let $x(t)$ be the solution of $\left\{
\begin{array}{c}
	\frac{dx(t)}{dt}=a(x(t),t) \\
	x(0)=x_{0}.
\end{array}
\right. $\\ 
Let $u=u(x(t),t)$,\\then the general solution of (9) can be obtained by using 
$$\frac{du}{dt}=b(x(t),t)u+\alpha(x(t),t)u^{n},$$ and taking $v=u^{1-n}$.\\
Let us consider the following examples:
\begin{Example} Let $a(x,t)=x$, $b(x,t)=t$, $\alpha(x,t)=x+t$ and $n=2$.
	\begin{equation}
		\left\{
		\begin{array}{c}
			u_{t}+xu_{x}=tu+(x+t)u^{2} \\
			u(x,0)=\phi(x).
		\end{array}
		\right.  
	\end{equation}
\end{Example}
We have $x(t)=x_{0}e^{t}$,
and $$\frac{du}{dt}=tu+(x(t)+t)u^{2}$$,
where $u=u(x(t),t)$.
Let $v=u^{-1}$, we get
$$\frac{dv}{dt}+tv=-(x(t)+t),$$
and 
$$\frac{d}{dt}(e^{\frac{t^{2}}{2}}v)=-(x(t)+t)e^{\frac{t^{2}}{2}},$$
$$v=-e^{-\frac{t^{2}}{2}}x_{0}\int_{0}^{t}e^{s+\frac{s^{2}}{2}}ds-e^{-\frac{t^{2}}{2}}(e^{\frac{t^{2}}{2}}-1)+
\frac{e^{-\frac{t^{2}}{2}}}{\phi(x_{0})}.$$
Then $$u(x,t)=(-e^{-\frac{t^{2}}{2}-t}x\int_{0}^{t}e^{s+\frac{s^{2}}{2}}ds+(e^{-\frac{t^{2}}{2}}-1)+
\frac{e^{-\frac{t^{2}}{2}}}{\phi(xe^{-t})})^{-1}.$$
\begin{figure}[hbt!]
	\centering
	\includegraphics[width=7cm]{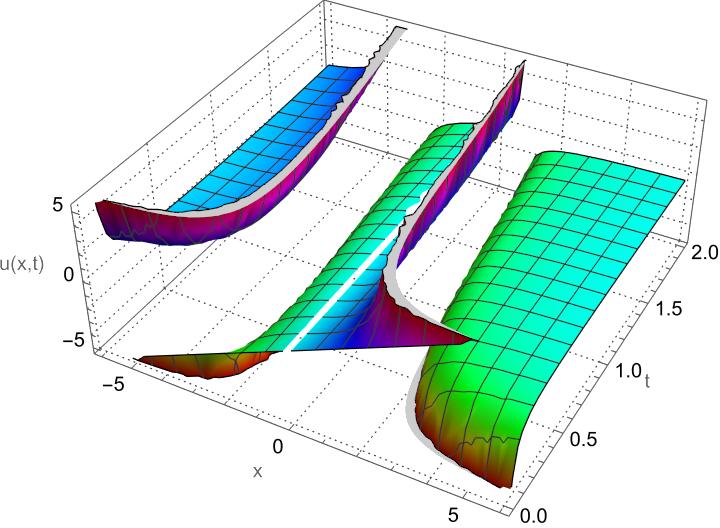} 
	\caption{Solution of (10) with $\phi(x)=x$.  }		
\end{figure}

\begin{Example} Let $a(x,t)=x$, $b(x,t)=1$, $\alpha(x,t)=t$ and $n=2$.
	\begin{equation}
		\left\{
		\begin{array}{c}
			u_{t}+xu_{x}=u+tu^{2} \\
			u(x,0)=\phi(x).
		\end{array}
		\right.  
	\end{equation}
\end{Example} 
Let $u=u(x(t),t)$,

We get $$\frac{du}{dt}-u=tu^{2}.$$
Let $v=u^{-1}$, then
$$\frac{dv}{dt}+v=-t,$$
and $$ \frac{d}{dt}(e^{t}v)=-t e^{t}.$$
Hence we obtain $$e^{t}v=\int_{0}^{t}-se^{s}ds+v(x_{0},0),$$  
$$v(x,t)=-t+1-e^{-t}+\frac{e^{-t}}{\phi(xe^{-t})},$$ and 
$$u(x,t)=(-t+1-e^{-t}+\frac{e^{-t}}{\phi(xe^{-t})})^{-1}.$$
\begin{figure}[hbt!]
	\centering
	\includegraphics[width=8cm]{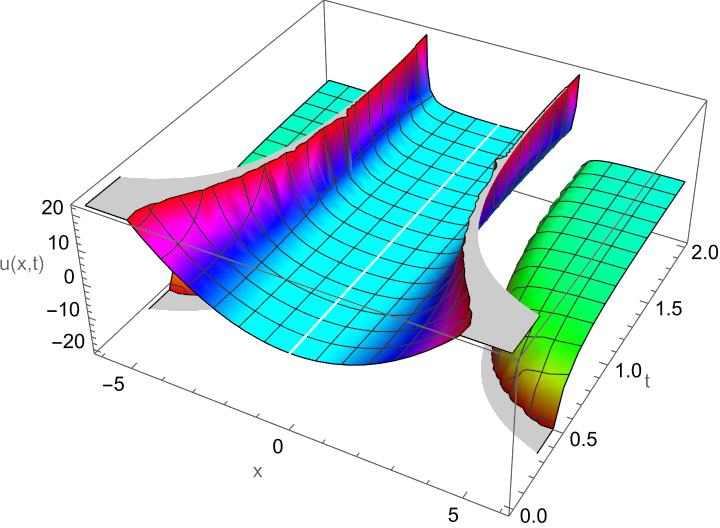} 
	\caption{Solution of (11)  with $\phi(x)=x^{2}.$  }	
\end{figure}
\begin{Example} Let $a(x,t)=1$, $b(x,t)=t$, $\alpha(x,t)=x$ and $n=2$.
	\begin{equation}
		\left\{
		\begin{array}{c}
			u_{t}+u_{x}=tu+xu^{2} \\
			u(x,0)=\phi(x).
		\end{array}
		\right. 
	\end{equation}
\end{Example}
Let $x(t)=t+x_{0}$ and $u=u(x(t),t),$
$$\frac{du}{dt}=tu+(t+x_{0})u^{2}.$$
Let $v=u^{-1}$, we get 
$$\frac{dv}{dt}+tv=-(t+x_{0}),$$
$$\frac{d}{dt}(e^{\frac{t^{2}}{2}}v)=-(t+x_{0})e^{\frac{t^{2}}{2}},$$ and 
$$e^{\frac{t^{2}}{2}}v=-\int_{0}^{t}(s+x_{0})e^{\frac{s^{2}}{2}}ds+\frac{1}{\phi(x-t)}.$$
Then $$u(x,t)=(-1+e^{-\frac{t^{2}}{2}}-(x-t)e^{-\frac{t^{2}}{2}}\int_{0}^{t}e^{\frac{s^{2}}{2}}ds+
\frac{e^{-\frac{t^{2}}{2}}}{\phi(x-t)})^{-1}.$$
\begin{figure}[hbt!]
	
	\centering
	\includegraphics[width=7cm]{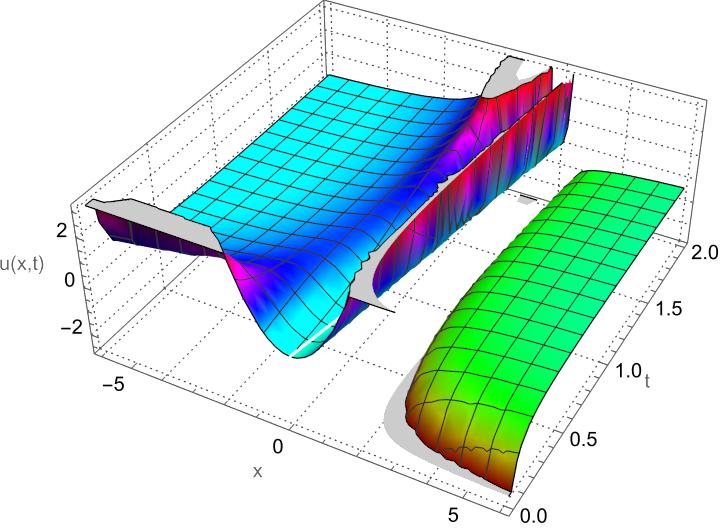} 
	
	\caption{Solution of (12)  with $\phi(x)=x^{2}.$  }	
\end{figure}

\section {Riccati Model	$
		u_{t}+a(x,t)u_{x}=b(x,t)u+\alpha(x,t)+\beta(x,t)u^{2} $ }
Let $u_{1}$ be the exact solution of Riccati Model 	$$u_{t}+a(x,t)u_{x}=b(x,t)u+\alpha(x,t)+\beta(x,t)u^{2}, $$then by using the substitution $u=u_{1}+w,$\\
we get the Bernoulli model 
$$w_{t}+a(x,t)w_{x}=\gamma(x,t)w+\beta(x,t)w^{2} $$ which is solved in the previous section.
\section{Partial Differential Equations of the form 
	$u_{tt}+a(x,t)u_{xt}=b(x,t)+(u_{t}+a(x,t)u_{x})f(u) $}
We consider the second order partial differential equation of the following type:
\begin{equation}	u_{tt}+a(x,t)u_{xt}=b(x,t)+(u_{t}+a(x,t)u_{x})f(u) \end{equation}
First let $x(t)$ be the solution of 
\begin{center}
	$\left\{
	\begin{array}{c}
		\frac{dx(t)}{dt}=a(x(t),t) \\
		x(0)=x_{0}.
	\end{array}
	\right. $ 
\end{center}
Then equation (13) takes the form
\begin{equation}\frac{d}{dt}(u_{t}(x(t),t))=b(x,t)+(u_{t}+a(x,t)u_{x})f(u).\end{equation}
We suppose that $u_{t}(x(t),t))=H(t)+K(u).$
\\ Then (14) is rewritten as $$\frac{d}{dt}(u_{t}(x(t),t))=H'(t)+K'(u)(u_{t}+au_{x}).$$
It is clear that (14) will take place if 
\begin{center}
	$\left\{
	\begin{array}{c}
		H'(t)=b(x,t)\\
		K'(u)=f(u).
	\end{array}
	\right. $ 
\end{center}
Hence, we have the following result:
\begin{Proposition}
	The second order partial differential equation (13) can be transformed to $$\frac{d}{dt}(u_{t}(x(t),t))=H'(t)+K'(u)(u_{t}+au_{x}),$$  where the functions $H$ and $K$ are the general solutions of $\left\{
	\begin{array}{c}
		H'(t)=b(x,t)\\
		K'(u)=f(u)
	\end{array}
	\right. $ 
\end{Proposition}
Let us consider the following example:
\begin{Example}  Let $a(x,t)=x$, $f(u)=u^{2}$ and $b(x,t)=x+t$.
	$$u_{tt}+xu_{xt}=x+t+(u_{t}+xu_{x})u^{2}. $$
\end{Example}
First, we solve 
	$\left\{
	\begin{array}{c}
		\frac{dx(t)}{dt}=x \\
		x(0)=x_{0}.
	\end{array}
	\right. $ 
	
The solution is $x=x_{0}e^{t}.$
The functions $H$ and $K$ are general solutions of
\begin{center}
	\begin{tabular}{ccc}
		$H'(t)$&=& $x(t)+t$,\\
		$K'(u)$&=& $u^{2}.$
	\end{tabular}
\end{center}
Then we get
$$H(t)=x_{0}e^{t}+\frac{t^{2}}{2}+C_{1},$$
$$H(t)=x+\frac{t^{2}}{2}+C_{1},$$
and $$K(u)=\frac{1}{3}u^{3}+C_{2}.$$
which reduces the PDEs to the form $$u_{t}(x,t)=x+\frac{t^{2}}{2}+\frac{u^{3}}{3}+C,$$
known as the Abel equation which can be solved by various methods. For more details, see ~\citep{D.Zwillinger,G.M.Murphy,Kamke}. In ~\citep{Panayotounakos},  the authors give implicit solutions using the first kind
Bessel's functions and the second kind of Newmann functions for the canonical form
of Abel's equation. 
\section{Partial Differential Equations of the form 
	$u_{xt}+a(x,t)u_{xx}=b(x,t)+(a(x,t)u_{x}+u_{t})f(u)$}
Let us consider the second order partial differential equation of the following type:
\begin{equation}u_{xt}+a(x,t)u_{xx}=b(x,t)+(a(x,t)u_{x}+u_{t})f(u) 
\end{equation}
We take $u_{x}(x(t),t)=H(t)+K(u)$,
where the functions $H$ and $K$ are the general solutions of
\begin{center}
	\begin{tabular}{c c c }
		$H'(t)$& =&$b(x(t),t)$ \\
		$K'(u)$& =&$f(u)$.
	\end{tabular}
\end{center}
\begin{Example}
	Consider the second order partial differential equation
	\begin{center}
		$u_{xt}+xu_{xx}=x+t+(u_{t}+xu_{x})u^{2}. $
	\end{center}
\end{Example}
In (15), we take $b(x,t)=x+t$, $a(x,t)=x$ and $f(u)=u^{2}$\\
The functions $H$ and $K$ are the general solutions of
\begin{center}
	\begin{tabular}{c c c }
		$H'(t)$& =&$x_{0}t+t$ \\
		$K'(u)$& =&$u^{2} $.
	\end{tabular}
\end{center}
Then $$u_{x}(x,t)=x+\frac{t^{2}}{2}+\frac{u^{3}}{3}+C.$$
The latter is Abel's equation which is integrable by various methods known in the literature. See ~\citep{D.Zwillinger}, ~\citep{G.M.Murphy} and  ~\citep{Kamke} for more details.
\section{Partial Differential Equations of the form	$f'(u_{t})(u_{tt}+au_{xt})=B(x,t)+A(u) (u_{t}+au_{x})$}
Consider the non linear second order differential equation of general form \begin{equation}
	f'(u_{t})(u_{tt}+au_{xt})=B(x,t)+A(u)(u_{t}+au_{x})
\end{equation}	 
Let $f(u_{t})$ of the form $H(t)+K(u)$,
then 
$$\frac{d}{dt}(f(u_{t}))=H'(t)+K'(u)(u_{t}+au_{x})$$
\begin{equation}f'(u_{t})(u_{tt}+au_{xt})=H'(t)+K'(u)(u_{t}+au_{x})
\end{equation}
Therefore the following statement holds
\begin{Proposition} The general solution of the non linear second order differential equation (16) is obtained by solving (17) where the functions $H$ and $K$ are the general solutions of	$H'(t)=B(x(t),t)$	and 	$K'(u)=A(u)$.
\end{Proposition}
\section{Partial Differential  Equation of the form	$u_{tt}+au_{xt}+b(u)u_{t}(u_{t}+au_{x})=\alpha(x,t)e^{-\int b(u)du} +G(u) (u_{t}+au_{x})$}
Consider the non linear second order partial differential equation of the general form 
\begin{equation}u_{tt}+au_{xt}+b(u)u_{t}(u_{t}+au_{x})=\alpha(x,t)e^{-\int b(u)du} +G(u) (u_{t}+au_{x})
\end{equation}
In ~\citep{Mahouton}, the author investigated special cases of (18) when $u$ and $\alpha$ are one variable functions.\\
Multiplying both sides of (18) by $e^{\int b(u)du}$, we get $$(u_{tt}+au_{xt})e^{\int b(u)du}+b(u)u_{t}(u_{t}+au_{x})e^{\int b(u)du}=\alpha(x,t)+G(u)(u_{t}+au_{x})e^{\int b(u)du}$$
Then $$\frac{d}{dt}[u_{t}(x(t),t)e^{\int b(u)du}]=\alpha(x,t)+G(u)(u_{t}+au_{x})e^{\int b(u)du}$$
Hence the nonlinear second order differential equation (18) is easily solved if we suppose that 
\begin{equation} u_{t}=(H(t)+K(t))e^{-\int b(u)du}\end{equation}
Then,
\begin{center}$H'(t)+K'(t)(u_{t}+au_{x})=\alpha(x,t)+G(u)(u_{t}+au_{x})e^{\int b(u)du}$\end{center} where the functions $H$ and $K$ are solutions of
$$H'(t)=\alpha(x(t),t)$$
$$K'(u)=G(u)e^{\int b(u)du}$$
Hence we obtain the following result:
\begin{Proposition} 
	The solution of the nonlinear second order differential equation (18) is obtained by solving (19) where the functions $H$ and $K$ are the general solutions of $H'(t)=\alpha(x(t),t)$ and $K'(u)=G(u)e^{\int b(u)du}$
\end{Proposition}
\begin{Example}
	Let us consider the nonlinear second order differential equation (18) with 
	\begin{center}$b(u)=-\frac{1}{u}$,$ G(u)=2 u^{2}$, $\alpha=x+t$ and $a=x$,
	\end{center}
Then  \begin{equation}\left\{\begin{array}{c}u_{tt}+xu_{xt}-\frac{u_{t}}{u}(u_{t}+xu_{x})=(x+t)u+2u^{2}(u_{t}+xu_{x})\\
			u(x,0)=x \\ u_{t}(x,0)=x^{2}+x^{3}\end{array}\right.
	\end{equation}
\end{Example}
First let $x(t)$ be the solution of 
\begin{center}
	$\left\{
	\begin{array}{c}
		\frac{dx(t)}{dt}=x \\
		x(0)=x_{0}.
	\end{array}\right.$ \\
\end{center}
Then 
$x(t)=x_{0}e^{t}.$

Taking $$H'(t)=\alpha(x(t),t)=x_{0}e^{t}+t,$$ and
$$K'(u)=2u^{2}(\frac{1}{u})=2u$$
we get $$H(t)=x_{0}e^{t}+\frac{t^{2}}{2}$$
and  $$K(u)=u^{2}.$$
We obtain
$$u_{t}=(x+\frac{t^{2}}{2}+u^{2}+c)u, $$
which can be reduced to the usual canonical form of Abel equation of the first kind  \begin{equation}w'_{t}(x,t)= w^{3}(x,t)+k(x,t).\end{equation}
The latter is integrable by various methods known in the literature. See ~\citep{D.Zwillinger}, ~\citep{G.M.Murphy} and  ~\citep{Kamke} 
for a good compilation of techniques developed to solve (21) for particular expressions of $k(x,t)$.
\\
If we take the initial condition $u(x,0)=x$ and $u_{t}(x,0)=x^{2}+x^{3}$, then $$u_{t}(x,t)=(x+\frac{t^{2}}{2}+u^{2})u. $$
\begin{figure}[hbt!]
	\centering 
	\includegraphics[width=7cm]{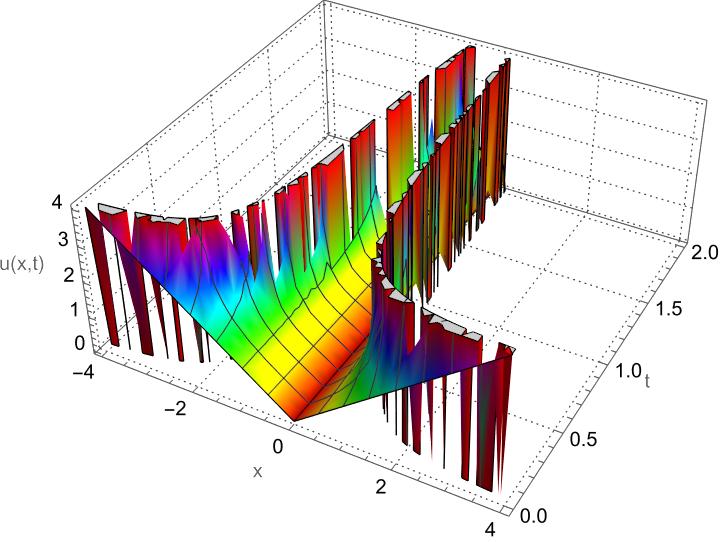} 
	\caption{Solution of (20)  with $u(x,0)=x$ and $u_{t}(x,0)=x^{2}+x^{3}.$  }	
\end{figure}
In the equation (18), we take
$b(u)=-\frac{1}{u}$ then
\begin{center}$u_{tt}+a(x,t)u_{xt}-\frac{u_{t}}{u}(u_{t}+a(x,t)u_{x})=\alpha(x,t)u+G(u)(u_{t}+a(x,t)u_{x}).$
\end{center}
If $G=u^{n}$,
\begin{equation}u_{tt}+a(x,t)u_{xt}-\frac{u_{t}}{u}(u_{t}+a(x,t)u_{x})=\alpha(x,t)u+u^{n}(u_{t}+a(x,t)u_{x}).\end{equation} 
Therefore we state the following result:
\begin{Proposition} 
	The exact solution of the non linear second order partial differential equation (22) is 
	$$u=[-e^{-n\int H(t)dt}\int e^{n\int H(t)dt} dt]^{-\frac{1}{n}}.$$
	where the function $H$ is the general solution of	$H'(t)=\alpha(x(t),t)$.
\end{Proposition}
Proof:
The equation (19) can be rewritten as \begin{equation}u_{t}=(H(t)+K(u))u \end{equation} 
where 
$H'(t)=\alpha(x(t),t)$ and 
$K'(u)=u^{n-1}.$

Substituting $K(u)=\frac{u^{n}}{n}$ to (23), we get
$$u_{t}=H(t)u+\frac{u^{n+1}}{n}$$ which is Bernoulli partial differential Equation where $H(t)$ is the general solution of $H'(t)=\alpha(x(t),t)$. Then we get
\begin{equation}u_{t}-H(t)u=\frac{u^{n+1}}{n}.\end{equation}
Let $v=u^{-n}$ then  $u=v^{-\frac{1}{n}}$ and 
$$u_{t}=-\frac{1}{n}v_{t}v^{-\frac{1}{n}-1}.$$
Equation (24) takes the form 
$$-\frac{1}{n}v_{t}v^{-\frac{1}{n}-1}-H(t)v^{-\frac{1}{n}}=\frac{1}{n}v^{-\frac{n+1}{n}},$$
which leads to the simpler form
$$v_{t}+nH(t)v=-1,$$
We obtain the solution   $$v=-e^{-n\int H(t)dt}\int e^{\int nH(t)dt} dt.$$
Finally, the exact solution of the nonlinear second order partial differential equation (22) is determined by $$u=[-e^{-n\int H(t)dt}\int e^{n\int H(t)dt} dt]^{-\frac{1}{n}}. $$

\section{Conclusion}
\vspace{0.3cm}
We have explored, in this paper, a large variety of partial differential equations which have been reduced to first order by variation of parameters combined with other techniques such as the method of characteristics. In most cases, the resulting first order differential equations are transformable to well-known integrable or solvable classical differential equations. We proved that it is possible to extend the techniques developed in this work to several classes of nonlinear second order partial differential equations. 

\bibliographystyle{unsrtnat}
\bibliography{references}  

\end{document}